\documentclass[final,leqno,onefignum,onetabnum]{siamltex1213}
\usepackage{enumitem}
\usepackage{lipsum}
\usepackage{graphicx}
\usepackage[makeroom]{cancel}

\usepackage{amsfonts}
\usepackage{amssymb,amsbsy}
\usepackage{amsmath}
\usepackage{amsfonts, amscd}
\usepackage{algorithm, algorithmic}
\usepackage{color}
\usepackage{url}
\usepackage{verbatim}
\usepackage{setspace}
\usepackage[abs]{overpic}
\usepackage{graphicx}
\usepackage{subfigure}
\usepackage{bbm, bm}
\usepackage{accents}
\usepackage{mleftright}
\usepackage{dsfont}

\newtheorem{example*}{Example}
\numberwithin{equation}{section}
\renewcommand{\theequation}{\thesection.\arabic{equation}}

\newtheorem{problem}{Problem}
\newtheorem{remark}{Remark}

\newtheorem{open question}{Open Question}

\newcommand{\mb}{\boldsymbol}
\renewcommand{\bm}{\boldsymbol}
\newcommand{\mc}{\mathcal}

\newcommand{\mcb}[1]{\mb{\mc{#1}}}
\newcommand{\norm}[2]{\left\| #1 \right\|_{#2}}
\newcommand{\reals}{\mathbb{R}}
\newcommand{\<}{\langle}
\renewcommand{\>}{\rangle}
\newcommand{\innerprod}[2]{\left\< #1, #2 \right\>}
\newcommand{\set}[1]{\left\{ #1 \right\}}
\newcommand{\eps}{\varepsilon}
\newcommand{\tensor}{\otimes}
\newcommand{\Tensor}{\bigotimes}

\DeclareMathOperator*{\argmin}{arg\,min}
\DeclareMathOperator*{\argmax}{arg\,max}

\newcommand{\wh}{\widehat}

\newcommand\Abs[1]{\mleft|#1\mright|}

\newcommand\citep\cite
\newcommand\citet\cite

\numberwithin{equation}{section}

\title{greedy approaches to \\
symmetric orthogonal tensor decomposition}

\author{Cun Mu\footnotemark[1]
  \and Daniel Hsu\footnotemark[2]
  \and Donald Goldfarb\footnotemark[1]}

\begin{document}
\maketitle
\slugger{simax}{xxxx}{xx}{x}{x--x}

{\renewcommand{\thefootnote}{\fnsymbol{footnote}}%
\footnotetext[1]{Department of Industrial Engineering and Operations Research, Columbia University (\url{cm3052@columbia.edu}, \url{goldfarb@columbia.edu}). DG was partially supported by NSF Grant CCF-1527809.}
\footnotetext[2]{Department of Computer Science, Columbia University (\url{djhsu@cs.columbia.edu}). DH was partially supported by NSF IIS-1563785, Bloomberg Data Science Research Grant and Sloan Research Fellowship.}
}

\begin{abstract}
Finding the symmetric and orthogonal decomposition (SOD) of a tensor is a recurring problem in signal processing, machine learning and statistics. In this paper, we review, establish and compare the perturbation bounds for two natural types of incremental rank-one approximation approaches. Numerical experiments and open questions are also presented and discussed.
\end{abstract}

\begin{keywords}
tensor decomposition, rank-1 tensor approximation, orthogonally decomposable tensor, perturbation analysis
\end{keywords}

\begin{AMS}
15A18, 15A69, 49M27, 62H25
\end{AMS}

\pagestyle{myheadings}
\thispagestyle{plain}
\markboth{C. MU, D. HSU, AND D. GOLDFARB}{Greedy Approaches to Tensor Decomposition}

\section{Introduction}\label{sec:intro}
A $p$-way $n$-dimensional tensor $\mcb T$, namely $\mcb T \in \bigotimes^p \reals^n := \reals^{n\times n \times \dotsb \times n}$, is called {\em symmetrically orthogonally decomposable (SOD)} \cite{mu2015successive, wang2016orthogonal} (a.k.a. {\em odeco} in \cite{robeva2016orthogonal}) if it can be expressed as a linear combination over the real field of symmetric $p$-th powers of $n$ vectors that generate an orthonormal basis of $\reals^n$. Mathematically, $\mcb T$ is SOD if there exist $\bm \lambda = [\lambda_1, \lambda_2, \ldots, \lambda_n] \in \reals^n$ and an orthogonal matrix $\mb V = [\bm v_1, \bm v_2, \ldots, \bm v_n] \in \reals^{n \times n}$ such that
\begin{flalign}\label{eqn:sodeco}
\mcb T = \lambda_1 \bm v_1^{\otimes p} + \lambda_2 \bm v_2^{\otimes p} + \cdots + \lambda_n \bm v_n^{\otimes p},
\end{flalign}
where $\bm v^{\otimes p}$, the symmetric $p$-th power of the vector $\bm v$, denotes a $p$-way $n$-dimensional tensor with $(\bm v^{\otimes p})_{i_1 i_2\cdots i_n} = v_{i_1}v_{i_2}\cdots v_{i_n}$. The decomposition $\set{(\lambda_i, \bm v_i)}_{i \in [n]}$ is called the {\em symmetric orthogonal decomposition (also abbreviated as SOD)} of $\mcb T$ with individual $\lambda_i$ and $\bm v_i$, respectively, called an {\em eigenvalue} and {\em eigenvector} of $\mcb T$.\footnote{For a more detailed discussion on eigenvalues and eigenvectors of SOD tensors, please see \cite{robeva2016orthogonal}.} The gist of our paper is to find the SOD of $\mcb T$ (potentially with perturbations), a recurring problem arising in different contexts including higher-order statistical
estimation~\cite{McCullagh87tensor}, independent component
analysis~\cite{comon1994independent, comon2010handbook}, and parameter
estimation for latent variable models \cite{JMLR:v15:anandkumar14b}, just to name a few.

From the expression \eqref{eqn:sodeco}, it is quite tempting to find $(\lambda_i, \bm v_i)$ one by one in a greedy manner using proper deflation procedures. Specifically, one first approximates $\mcb T$ by the best rank-one tensor,
\begin{flalign}\label{eqn:TBRO}
  (\lambda^\star, \bm v^\star) \in \argmin_{\lambda \in \reals, \norm{\bm v}{}
  = 1} \norm{\mcb T - \lambda \cdot \bm v^{\tensor p}}{F}.
\end{flalign}
After that, to find the next pair, one modifies the optimization problem \eqref{eqn:TBRO} to exclude the found eigenpair $(\lambda^\star, \bm v^\star)$.
We next review two natural deflation procedures---{\em residual deflation} \cite{zhang2001rank} and {\em constrained deflation} \cite{kolda2001orthogonal}---which incorporate the information of $(\lambda^\star, \bm v^\star)$ into an optimization framework by altering, respectively, the objective and the feasible set of problem \eqref{eqn:TBRO}.

\paragraph{Residual deflation} In residual deflation, the rank-one approximation is subtracted from the original tensor, i.e., $\mcb T \gets \mcb T - \lambda^\star \cdot (\bm v^\star)^{\tensor p}$, and then finds the best rank-one approximation to the deflated tensor by solving \eqref{eqn:TBRO} again. The complete scheme, referred to as  {\em Successive Rank-One Approximation with Residual Deflation (SROAwRD)}, is described in Algorithm \ref{alg:SROAwRD}.

\begin{algorithm}
\caption{\underline{S}uccessive \underline{R}ank-\underline{O}ne \underline{A}pproximation \underline{w}ith \underline{R}esidual \underline{D}eflation (SROAwRD)}
\label{alg:SROAwRD}
\begin{algorithmic}[1]
  \renewcommand\algorithmicrequire{\textbf{input}}

  \REQUIRE a symmetric $p$-way tensor $\wh{\mcb T} \in \bigotimes^p \reals^n$.
  \STATE \textbf{initialize} $\wh{\mcb T}_0 \gets \wh{\mcb T}$
  \FOR{$k=1$ to $n$}
    \STATE $(\hat\lambda_k,\hat{\bm v}_k) \in
    \argmin_{\lambda\in\reals,\norm{\bm v}{}=1} \; \norm{\wh{\mcb
    T}_{k-1} - \lambda \bm v^{\tensor p}}{F}$.
    \label{line:rank_one_approx_}

    \STATE $\wh{\mcb T}_k \gets \wh{\mcb T}_{k-1} - \hat{\lambda}_k
    \hat{\bm v}_k^{\tensor p}$.

  \ENDFOR

  \RETURN $\{ (\hat\lambda_k, \hat{\bm v}_k) \}_{k =1}^{n}$.
\end{algorithmic}
\end{algorithm}

\paragraph{Constrained deflation} In constrained deflation, one restricts $\bm v$ to be nearly orthogonal to $\pm \bm v^\star$ by solving problem \eqref{eqn:TBRO} with the additional linear constraints $-\theta \le \innerprod{\bm v^\star}{\bm v} \le \theta$, where $\theta > 0$ is a prescribed parameter. The complete scheme, referred to as {\em Successive Rank-One Approximation with Constrained Deflation (SROAwCD)}, is described in Algorithm \ref{alg:SROAwCD}. At the $k$-th iteration, rather than deflating the original tensor $\wh{\mcb T}$ by subtracting from it the sum of the $(k-1)$ rank-one tensors $\hat \lambda_1 \hat{\bm v}_1^{\tensor p}$, $\hat \lambda_2 \hat{\bm v}_2^{\tensor p}$, $\cdots$, $\hat \lambda_{k-1} \hat{\bm v}_{k-1}^{\tensor p}$ as the SROAwRD method does, the SROAwCD method imposes the near-orthogonality constraints $|\innerprod{\bm v}{\hat{\bm v}_{i}}| \le \theta$ for $i = 1,2,\ldots, k-1$.

\begin{algorithm}[ht]
\caption{\underline{S}uccessive \underline{R}ank-\underline{O}ne \underline{A}pproximation \underline{w}ith \underline{C}onstrained \underline{D}eflation (SROAwCD)}
\label{alg:SROAwCD}
\begin{algorithmic}[1]
  \renewcommand\algorithmicrequire{\textbf{input}}

  \REQUIRE a symmetric $p$-way tensor $\wh{\mcb T} \in \bigotimes^p \reals^n$, parameter $\theta > 0$.

  \STATE \textbf{initialize} $\hat{\bm v}_0 \gets \bm 0$
  \FOR{$k=1$ to $n$}
  \STATE Solve the following optimization problem:
           \begin{flalign} \label{eqn:SROAwCD}
           (\hat\lambda_k,\hat{\bm v}_k) \in &\argmin_{\lambda\in \reals, \bm v \in \reals^n} \quad \norm{\wh{\mcb T} - \lambda \bm v^{\tensor p}}{F} \\
           &\;\;\qquad\mbox{s.t.} \quad \norm{\bm v}{} = 1  \nonumber \\
           & \qquad \qquad \;\;\;\; -\theta \le \innerprod{\bm v}{\hat{\bm v}_i} \le \theta, \quad i = 0, \;1,\;2,\;\dots, \;k-1  \nonumber
           \end{flalign}
  \ENDFOR

  \RETURN $\{ (\hat\lambda_k, \hat{\bm v}_k) \}_{k = 1}^n$.
\end{algorithmic}
\end{algorithm}

It is not hard to prove that given the SOD tensor $\mcb T = \sum_{i \in [n]} \lambda_i \bm v_i^{\tensor p}$ as the input, both SROAwRD and SROAwCD methods are capable of finding the eigenpairs $\set{(\lambda_i, \bm v_i)}_{i \in [n]}$ exactly. In this paper, we focus on the more challenging case of tensors that are only close to being SOD. Specifically:
\begin{problem}\label{question:main}
Suppose the SOD tensor $\mcb T = \sum_{i \in [n]} \lambda_i \bm v_i^{\tensor p}$, and that the perturbed SOD tensor $\wh{\mcb T}$ is provided as input to the SROAwRD and SROAwCD methods. Characterize the discrepancy between $\{(\lambda_i, \bm v_i)\}_{i \in [n]}$ and the components $\{(\hat\lambda_i, \hat{\bm v}_i)\}_{i \in [n]}$ found by these methods.
\end{problem}

In this paper, we provide positive answers to Problem \ref{question:main}. The characterization for SROAwRD was done in our previous paper \cite{mu2015successive}; we review the results in Section \ref{sec:SROAwRD}. The charaterization for SROAwCD is the main contribution of the present paper. These results can be regarded as higher order generalizations of the Davis-Kahan perturbation result \cite{davis1970rotation} for matrix eigen-eigenvalue decomposition, and is not only of mathematical interest but also crucial to applications in signal processing, machine learning and statistics \cite{McCullagh87tensor, comon1994independent, comon2010handbook, JMLR:v15:anandkumar14b}, where the common interest is to find the underlying eigenpairs $\set{(\lambda_i, \bm v_i)}$ but the tensor collected is subject to inevitable perturbations arising from sampling errors, noisy measurements, model specification, numerical errors and so on.

\paragraph{Organization} The rest of the paper is organized as follows. In Section 2, we introduce notation relevant to this paper. In Section 3, we review theoretically what is known about the SROAwRD method. In Section 4, we provide a perturbation analysis for the SROAwCD method.

\section{Notation}\label{sec: notation}
In this section, we introduce some tensor notation needed in our paper, largely borrowed from \cite{lim2005eigenvalue}.

\paragraph{Symmetric tensor} A real $p$-way $n$-dimensional tensor
$\mcb A \in \Tensor^p \reals^n := \reals^{n\times n \times \dotsb \times n}$,
\[
\mcb A = \left(\mc A_{i_1, i_2, \dotsc, i_p} \right), \;\; \mc A_{i_1,
i_2, \dotsc, i_p} \in \reals, \quad 1\le i_1, i_2, \dotsc, i_p \le n,
\]
is called \emph{symmetric} if its entries are invariant under any
permutation of their indices, i.e. for any $i_1, i_2, \dotsc, i_p \in [n]:=\set{1,2,\ldots n}$,
\[
  \mc A_{i_1, i_2, \dotsc, i_p} = \mc A_{i_{\pi(1)}, i_{\pi(2)},
  \dotsc, i_{\pi(p)}}
\]
for every permutation mapping $\pi$ of $[p]$.
\paragraph{Multilinear map} In addition to being considered as a multi-way array, a
tensor $\mcb A \in \Tensor^p \reals^n$ can also be interpreted as a
multilinear map in the following sense: for any matrices $\mb V_i \in
\reals^{n \times m_i}$ for $i \in [p]$, we define
$\mcb A(\mb V_1, \mb V_2, \dotsc, \mb V_p)$ as a tensor in
$\reals^{m_1\times m_2 \times \cdots \times m_p}$ whose $(i_1, i_2,
\dotsc, i_p)$-th entry is
\[
\left(\mcb A(\mb V_1, \mb V_2, \dotsc, \mb V_p)\right)_{i_1, i_2,
\dotsc, i_p} := \sum_{j_1, j_2, \dotsc, j_p \in [n]} \mc A_{j_1, j_2,
\dotsc, j_p} (V_1)_{j_1 i_1} (V_2)_{j_2 i_2} \cdots  (V_p)_{j_p i_p}.
\]
The following two special cases are quite frequently used in the paper:
\vspace{5mm}
\begin{description}[leftmargin=12.5pt]
  \item[{$\rhd$}] $\mb V_i = \bm x \in \reals^n$ for all $i \in [p]$: $\mcb A \bm x^{\tensor p} := \mcb A(\bm x, \bm x, \dotsc, \bm x),$
    which defines a  homogeneous polynomial of degree $p$.
  \vspace{2mm}
  \item[{$\rhd$}] $\mb V_i = \bm x \in \reals^n$ for all $i \in [p-1]$, and $\mb
    V_p = \mb I \in \reals^{n \times n}$:
    \begin{align*}
      \mcb A \bm x^{\tensor p-1}
      & :=
      \mcb A(\bm x,\dotsc,\bm x, \mb I)
      \in \reals^n.
    \end{align*}

\end{description}
For a symmetric tensor $\mcb A \in \Tensor^p \reals^n$,
     the differentiation result $\nabla_{\bm x} \left(\mcb A \bm x^{\tensor p}\right) = p\cdot \left(\mcb A \bm x^{\tensor p-1}\right)$ can be established.

\paragraph{Inner product}
For any tensors $\mcb A$, $\mcb B \in \Tensor^p \reals^n$, the inner
product between them is naturally defined as
\[
  \innerprod{\mcb A}{\mcb B} := \sum_{i_1,i_2, \dotsc,i_p \in [n]}
  \mc A_{i_1, i_2, \dotsc, i_p}  \mc B_{i_1, i_2, \dotsc, i_p}.
\]
\paragraph{Tensor norms}
Two tensor norms will be used in the paper.
For a tensor $\mcb A \in \Tensor^p \reals^n$, its \emph{Frobenius
norm} is $\norm{\mcb A}{F} := \sqrt{\innerprod{\mcb A}{\mcb A}}$, and
its \emph{operator norm} $\norm{\mcb A}{}$, is defined as
$\max_{\norm{\bm x_i}{} = 1} \mcb A(\bm x_1, \bm x_2, \dotsc, \bm
x_p)$.
It is also well-known that for symmetric tensors $\mcb A$,
$\norm{\mcb A}{} $ can be equivalently defined as $ \max_{\norm{\bm x}{}=1} |\mcb A
\bm x^{\tensor p}|$ (see, e.g.,~\cite{chen2012maximum,
zhang2012best}).

\section{Review on SROAwRD}\label{sec:SROAwRD}
Algorithm \ref{alg:SROAwRD} is intensively studied in the tensor community, though most papers \cite{de2000best, zhang2001rank,
kofidis2002best, wang2007successive, kolda2011shifted,
han2012unconstrained, chen2012maximum, zhang2012best, l2014sequential, jiang2012tensor, nie2013semidefinite,
yang2014properties, JMLR:v15:anandkumar14b, hu2016note, 7473918} focus on the numerical aspects of how to solve the best tensor rank-one approximation \eqref{eqn:TBRO}. Regarding theoretical guarantees for the symmetric and orthogonal decomposition, Zhang and Golub \cite{zhang2001rank} first prove that SROAwRD outputs the exact symmetric and orthogonal decomposition if the input tensor is  symmetric and orthogonally decomposable:

\begin{proposition}{\cite[Theorem 3.2]{zhang2001rank}}\label{thm:main-SROA}
Let $\mcb T \in \Tensor^p \reals^n$ be a symmetric tensor with orthogonal decomposition $\mcb T = \sum_{i\in [n]} \lambda_i \bm v_i^{\tensor p}$, where $\lambda_i \neq 0$ and  $\set{\bm v_1, \bm v_2 ,\ldots , \bm v_n}$ forms an orthonormal basis of $\reals^n$. Let $\{ (\hat\lambda_i, \hat{\bm v}_i) \}_{i \in [n]}$ be the output of Algorithm \ref{alg:SROAwRD} with input $\mcb T$. Then $\mcb T = \sum_{i\in [n]} \hat\lambda_i \hat{\bm v}_i^{\tensor p}$, and moreover there exists a permutation $\pi$ of $[n]$ such that for each $j\in [n]$,
\begin{flalign*}
&\min\; \set{|\lambda_{\pi(j)} - \hat{\lambda}_j|, \;|\lambda_{\pi(j)} + \hat{\lambda}_j|} = 0, \\
&\min\; \set{\norm{\bm v_{\pi(j)} - \hat{\bm v}_j}{},\; \norm{\bm v_{\pi(j)} + \hat{\bm v}_j}{}} = 0.
\end{flalign*}
\end{proposition}

The perturbation analysis is recently addressed in \cite{mu2015successive}:

\begin{theorem}{\cite[Theorem 3.1]{mu2015successive}}\label{thm:main-SROA-E}
  There exists a positive constant $c$ such that the
  following holds.
  Let $\wh{\mcb T} := \mcb T + \mcb E \in \Tensor^p \reals^n$, where the ground truth tensor
  $\mcb T$ is symmetric with orthogonal decomposition $\mcb T = \sum_{i \in [n]} \lambda_i \bm v_i^{\tensor p}$, $\set{\bm v_1, \bm v_2 ,\ldots , \bm v_n}$ forms an orthonormal basis of $\reals^n$, $\lambda_i \neq 0$ and the perturbation tensor $\mcb E$ is symmetric
  with operator norm $\eps := \norm{\mcb E}{}$.
  Assume $\eps \le c \cdot \lambda_{\min}/n^{1/(p-1)}$, where
  $\lambda_{\min} := \min_{i\in[n]} |\lambda_i|$.
  Let $\{ (\hat\lambda_i,\hat{\bm v}_i) \}_{i \in [n]}$ be the output
  of Algorithm~\ref{alg:SROAwRD} with input $\wh{\mcb T}$.
  Then there exists a permutation $\pi$ over $[n]$ such that for each $j\in [n]$:
\begin{flalign*}
&\min\; \set{|\lambda_{\pi(j)} - \hat{\lambda}_j|, \;|\lambda_{\pi(j)} + \hat{\lambda}_j|} \le 2 \eps, \\
&\min\; \set{\norm{\bm v_{\pi(j)} - \hat{\bm v}_j}{},\; \norm{\bm v_{\pi(j)} + \hat{\bm v}_j}{}} \le 20 \eps / \Abs{\lambda_{\pi(j)}}.
\end{flalign*}
\end{theorem}

Theorem \ref{thm:main-SROA-E} generalizes Proposition \ref{thm:main-SROA}, and provides perturbation bounds for the SROAwRD method. Specifically, when the operator norm of the perturbation tensor vanishes, i.e. $\eps = 0$, Theorem \ref{thm:main-SROA-E} is reduced to Proposition \ref{thm:main-SROA}; when $\eps$ is small enough (i.e. $\eps = O(1/n^{1/(p-1)})$), the SROAwRD method is able to robustly recover the eigenpairs $\set{(\lambda_i, \bm v_i)}_{i\in[n]}$ of the underlying symmetric and orthogonal decomposable tensor $\mcb T$.

In Theorem \ref{thm:main-SROA-E}, $\eps$ is required to be at most on the order of $1/n^{1/(p-1)}$, which decreases with increasing tensor size. It is interesting to explore whether or not this dimensional dependency is essential:
\begin{open question}
Can we provide a better analysis for the SROAwRD method to remove the dimensional dependance on the noise level? Or can we design a concrete example to corroborate the necessity of this dimensional dependency?
\end{open question}

The existence of the dimensional dependency, at least for the current proof in Mu et. al. \cite{mu2015successive}, can be briefly explained as follows.
At the end of the $k$-th iteration, we subtract the rank-one tensor $\hat \lambda_k \hat{\bm v}_k^{\tensor p}$ from $\wh{\mcb T}_{k-1}$. Since $(\hat \lambda_k, \hat{\bm v}_k)$  only approximates the underlying truth, this deflation procedure introduces additional errors into $\wh{\mcb T}_k$. Although \cite{mu2015successive} has made substantial efforts to reduce the accumulative effect from sequential deflation steps, the perturbation error $\eps$ still needs to depend on the iteration number in order to control the perturbation bounds of the eigenvalue and eigenvector, and we tend to believe that the dimensional dependency in Theorem \ref{thm:main-SROA-E} is necessary.

In contrast, the SROAwCD method, instead of changing the objective, imposes additional constraints, which force the next eigenvector $\hat{\bm v}_{k}$ to be
nearly orthogonal to $\set{\hat{\bm v}_1, \hat{\bm v}_2,\ldots, \hat{\bm v}_{k-1}}$.
As the SROAwCD method alters the search space rather than the objective in the optimization, there is hope that the requirement on the noise level might be dimension-free. In the next section, we will confirm this intuition.

\section{SROAwCD}\label{sec:SROAwCD}
In this section, we establish the first perturbation bounds that have been given for the SROAwCD method for tensor SOD. The main result can be stated as follows:
\begin{theorem}\label{thm:C-SROA}
Let $\wh{\mcb T} := \mcb T + \mcb E \in \Tensor^p \reals^n$, where $\mcb T$ is a symmetric tensor with orthogonal decomposition $\mcb T
  = \sum_{i=1}^n \lambda_i \bm v_i^{\tensor p}$, $\set{\bm v_1, \bm
  v_2, \dotsc, \bm v_n}$ is an orthonormal basis of $\reals^n$,
  $\lambda_i \neq 0$ for all $i \in [n]$, and $\mcb E$ is a symmetric
  tensor with operator norm $\eps:=\norm{\mcb E}{}$.
  Assume $0 < \theta\le {1}/{(2 \kappa)}$ and $\eps\le \theta^2 \lambda_{\min}/ 12.5$, where  $\kappa := \lambda_{\max}/\lambda_{\min}$, $\lambda_{\min} := \min_{i\in[n]} \Abs{\lambda_i}$ and $\lambda_{\max} := \max_{i\in[n]} \Abs{\lambda_i}$. Let $\{ (\hat\lambda_i,\hat{\bm v}_i) \}_{i \in [n]}$ be the output of Algorithm \ref{alg:SROAwCD} for input $(\wh{\mcb T},\theta)$. Then there exists a permutation $\pi$ of $[n]$ such that for all $j \in [n]$,
\begin{flalign}
&\min\; \set{|\lambda_{\pi(j)} - \hat{\lambda}_j|, \;|\lambda_{\pi(j)} + \hat{\lambda}_j|} \le \eps, \label{eqn:lam_bd}\\
&\min\; \set{\norm{\bm v_{\pi(j)} - \hat{\bm v}_j}{},\; \norm{\bm v_{\pi(j)} + \hat{\bm v}_j}{}} \le (6.2 + 4 \kappa) \eps/ |\lambda_{\pi(j)}|. \label{eqn:v_bd}
\end{flalign}
\end{theorem}

Theorem \ref{thm:C-SROA} guarantees that for an appropriately chosen $\theta$, the SROAwCD method can approximately recover $\set{(\lambda_i, \bm v_i)}_{i \in [n]}$ whenever the perturbation error $\eps$ is small. A few remarks immediately come to find. First, Theorem \ref{thm:C-SROA} specifies the choice of the parameter $\theta$, which depends on the ratio of the largest to smallest eigenvalues of $\mcb T$ in absolute value. In subsection 4.2, we will see this dependency is necessary through numerical studies. Second, in contrast to the SROAwRD method, Theorem \ref{thm:C-SROA} does not require the noise level to be dependent on the tensor size. This could be a potential advantage for the SROAwCD method.

The rest of this section is organized as follows. In subsection \ref{subsec:proof}, we provide the proof for Theorem \ref{thm:C-SROA}. In subsection \ref{subsec:numeric}, we present numerical experiments to corroborate Theorem \ref{thm:C-SROA}. In subsection \ref{subsec:disc}, we discuss issues related to determining the maximum spectral ratio $\kappa$ defined in Theorem \ref{thm:C-SROA}.

\subsection{Proof of Theorem \ref{thm:C-SROA}} \label{subsec:proof}
We will prove Theorem \ref{thm:C-SROA} by induction. For the base case, we need the perturbation result regarding the best rank-one tensor approximation, which is proven in \cite{mu2015successive} and can be regarded as a generalization of its matrix counterpart \cite{weyl1912asymptotische, davis1970rotation}. In the following, we restate this result \cite[Theorem 2.2]{mu2015successive} with a minor variation:
\begin{lemma}\label{lem:for_base}
Let $\wh{\mcb T} := \mcb T + \mcb E \in \Tensor^p \reals^n$, where $\mcb T$ is a symmetric tensor with orthogonal decomposition $\mcb T
  = \sum_{i=1}^n \lambda_i \bm v_i^{\tensor p}$, $\set{\bm v_1, \bm
  v_2, \dotsc, \bm v_n}$ is an orthonormal basis of $\reals^n$,
  $\lambda_i \neq 0$ for all $i \in [n]$, and $\mcb E$ is a symmetric
  tensor with operator norm $\eps:=\norm{\mcb E}{}$. Let $(\hat \lambda, \hat{\bm v}) \in \argmin_{\lambda\in \reals, \norm{\bm v}{} = 1}\; \norm{\wh{\mcb T} - \lambda \bm v^{\tensor p}}{F}$. Then there exist $j \in [n]$ such that
  \begin{flalign*}
&\min\; \set{|\lambda_{j} - \hat{\lambda}|, \;|\lambda_{j} + \hat{\lambda}|} \le \eps, \quad \mbox{and} \\
&\min\; \set{\norm{\bm v_{j} - \hat{\bm v}}{},\; \norm{\bm v_{j} + \hat{\bm v}}{}} \le 10 \left(\frac{\eps}{|\lambda_j|} +  \left(\frac{\eps}{\lambda_j}\right)^2 \right).
\end{flalign*}
\end{lemma}

Now we are ready to prove our main Theorem \ref{thm:C-SROA}.

\begin{proof}
Without loss of generality, we assume $p\ge 3$ is odd, and $\lambda_i > 0$ for all $i \in [n]$ (as we can always flip the signs of the $\bm v_i's$ to ensure this). Then problem \eqref{eqn:SROAwCD} can be equivalently written as
\begin{flalign}\label{eqn:equiv_odd}
\hat{\bm v}_k \in \argmax_{\bm v \in \reals^n} \;\wh{\mcb T} \bm v^{\tensor p} \qquad \mbox{s.t.} \quad \norm{\bm v}{} = 1, \; \mbox{and} \; |\innerprod{\bm v}{\hat{\bm v}_i}| \le \theta \;\; \forall \; i \in [k-1],
\end{flalign}
and $\hat \lambda_k = \wh{\mcb T} \hat{\bm v}_k^{\tensor p}$.

To prove the theorem, it suffices to prove that the following property holds for each $k \in [n]$: there is a permutation $\pi$ of $[n]$ such that for every $j \in [k]$,
  \begin{flalign}
      |\lambda_{\pi(j)} - \hat{\lambda}_j|
      \le \eps \quad \text{and} \quad
      \norm{\bm v_{\pi(j)}-\hat{\bm v}_j}{}
      \le \frac{(6.2 + 4\kappa) \eps}{\lambda_{\pi(j)}}.
      \tag{$*$}
      \label{eqn:indhyp}
  \end{flalign}
We will prove \eqref{eqn:indhyp} by induction.

For the base case $k = 1$, Lemma \ref{lem:for_base} implies that
there exists a $j \in [n]$ satisfying
\[
|\hat \lambda_1 - \lambda_j| \le \eps, \quad \mbox{and} \quad \norm{\hat{\bm v}_1 - \bm v_j}{} \le 10\frac{\eps}{\lambda_j}\left(1+ \frac{\eps}{\lambda_j} \right)\le \frac{10.2 \eps}{\lambda_j} \le (6.2 + 4\kappa) \frac{\eps}{\lambda_j},
\]
where we have used the fact that $\eps/\lambda_j \le \eps/\lambda_{\min} \le \theta^2/12.5 \le 1/50.$

Next we assume the induction hypothesis \eqref{eqn:indhyp} is true for $k \in [n-1]$, and prove that there exists an $l \in [n] \backslash \set{\pi(j): j\in [k]}$ that satisfies
\begin{flalign}\label{eqn:indhyp_l}
|\hat \lambda_{k+1} - {\lambda}_l|\le \eps, \quad \mbox{and} \quad \norm{\hat{\bm v}_{k+1}- {\bm v}_l}{}
      \le \frac{(6.2 + 4\kappa) \eps}{\lambda_{\pi(l)}}.
\end{flalign}

Denote $\hat{\bm x}:= \hat{\bm v}_{k+1}$ and $\hat \lambda := \hat \lambda_{k+1}$. Then based on \eqref{eqn:equiv_odd}, one has
\begin{flalign}\label{eqn:subproblem}
 \hat{\bm x} \;\in\; \arg \max_{ \bm v \in \reals^n} \;\; \wh{\mcb T} {\bm v}^{\tensor p}  \quad \mbox{s.t.} \;\; \norm{\bm v}{} = 1, \; |\innerprod{\hat{\bm v}_i}{\bm v}| \le \theta\;\; \forall  i \in [k],
\end{flalign}
and $\hat \lambda = \wh{\mcb T} \hat{\bm x}^{\tensor p}$. Since $\set{\bm v_i}_{i \in [n]}$ forms an orthonormal basis, we may write $\hat{\bm x} = \sum_{i \in [n]} x_i \bm v_i$. Without loss of generality, we renumber $\set{\left(\lambda_{\pi(i)}, \bm v_{\pi(i)} \right)}_{i \in [k]}$ to $\set{\left(\lambda_i, \bm v_i \right)}_{i \in [k]}$ and renumber $\set{\left(\lambda_i, \bm v_i \right)}_{i \in [n] \backslash \set{\pi(i) \vert i \in [k]}}$ to $\set{\left(\lambda_i, \bm v_i \right)}_{i \in [n]\backslash[k]}$, respectively, to satisfy
\begin{flalign}\label{eqn:ordering_ass}
&\lambda_1 |x_1|^{p-2} \ge \lambda_2 |x_2|^{p-2} \ge \ldots \ge \lambda_k |x_k|^{p-2}, \quad \mbox{and}  \\
&\lambda_{k+1}|x_{k+1}|^{p-2} \ge \lambda_{k+2}|x_{k+2}|^{p-2} \ge \ldots \ge \lambda_{n}|x_{n}|^{p-2}. \nonumber
\end{flalign}

In the following, we will show that $l = k+1$ is indeed the index satisfying \eqref{eqn:indhyp_l}. The idea of the rest of the proof is as follows. We first provide lower and upper bounds for $\hat \lambda = \wh{\mcb T} \hat{\bm x}^{\tensor p}$, based on which, we are able to show that $|\hat{\lambda} - \lambda_{l}| = O(\eps)$ and $1 - |\innerprod{\hat{\bm v}}{\bm v_{l}}| = O(\eps)$. However, Theorem \ref{thm:C-SROA} requires $1 - |\innerprod{\hat{\bm v}}{\bm v_{l}}| = O(\eps^2)$. To close this gap, we characterize the optimality condition of \eqref{eqn:subproblem}, use of which enables us to sharpen the upper bound of $1 - |\innerprod{\hat{\bm v}}{\bm v_l}|$.

We first consider the lower bound for $\hat \lambda$ by finding a $\bm v$ that is feasible for \eqref{eqn:subproblem}. For each $(i,j) \in [n]\backslash [k] \times [k]$, one has
\begin{flalign} \label{eqn:v_feas}
|\innerprod{\bm v_i}{\hat{\bm v}_{j}}| & = |\innerprod{\bm v_i}{\hat{\bm v}_j - \bm v_{\pi(j)}}| \le \norm{\hat{\bm v}_j - \bm v_{\pi(j)}}{} \\
& \le  \frac{(6.2 + 4 \kappa)\eps}{\lambda_{\pi(j)}} \le \frac{(6.2 + 4 \kappa) \theta^2 \cancel{\lambda_{\min}}}{12.5 \cancel{\lambda_{\min}}} = \frac{(6.2 + 4 \kappa) \theta}{12.5}\cdot \theta = \frac{6.2 + 4\kappa}{25\kappa}\theta < {\theta}, \nonumber
\end{flalign}
where we have used the Cauchy-Schwarz inequality, and the facts $\pi(j) \in [k]$, $\innerprod{\bm v_i}{\bm v_{\pi (j)}} = 0$, $\eps \le {\theta^2 \lambda_{\min}}/{12.5}$ and $\theta\le {1}/{(2\kappa)}$.
Hence, $\set{\bm v_i}_{i \in [n]\backslash [k]}$ are all feasible to problem \eqref{eqn:subproblem} and then we can naturally achieve a lower bound for $\hat \lambda$, as
\begin{flalign}\label{eqn:lower_bd}
\hat \lambda = \wh{\mcb T} \hat{\bm x}^{\tensor p} \ge \max_{i \in [n]\backslash [k]} \wh{\mcb T} \bm v_i^{\tensor p} \ge \max_{i \in [n]\backslash [k]} \lambda_i - \eps \ge \lambda_{k+1} - \eps.
\end{flalign}

Regarding the upper bound for $\hat \lambda$, one has
\begin{flalign}
\hat \lambda = \wh{\mcb T} \hat{\bm x}^{\tensor p} = (\mcb T + \mcb E) \hat{\bm x}^{\tensor p} & =  \left(\sum_{i=1}^n \lambda_i \bm v_i^{\tensor p} + \mcb E \right) \hat{\bm x}^{\tensor p}  \nonumber\\
&= \sum_{i = 1}^k \lambda_i x_i^p + \sum_{i = k+1}^n \lambda_i x_i^p   + \mcb E \hat{\bm x}^{\tensor p} \nonumber \\
& \le \sum_{i = 1}^k \lambda_i |x_i|^{p-2} x_i^2 + \sum_{i = k+1}^n \lambda_i |x_i|^{p-2}x_i^2 + \mcb E \hat{\bm x}^{\tensor p} \nonumber \\
& \le \max\{\lambda_1 |x_1|^{p-2}, \lambda_{k+1} |x_{k+1}|^{p-2}\} + \eps, \label{eqn:upper_bd}
\end{flalign}
where the last line is due to the assumptions made in \eqref{eqn:ordering_ass}, $\norm{\bm x}{} = 1$ and $\eps = \norm{\mcb E}{}$.

Combining \eqref{eqn:lower_bd} and \eqref{eqn:upper_bd}, we have
\begin{flalign}\label{eqn:comb_upp_low}
\lambda_{k+1} - \eps \le \max_{i \in [n]\backslash [k]} \lambda_i - \eps  \le \hat \lambda \le \max\{\lambda_1 |x_1|^{p-2}, \lambda_{k+1} |x_{k+1}|^{p-2}\} + \eps.
\end{flalign}
Also note that
\begin{flalign}
&\lambda_1 |x_1|^{p-2} + \eps \\
&\le  \lambda_1 |x_1| + \eps
= \lambda_1 |\innerprod{\hat{\bm x}}{\bm v_1}| + \eps
= \lambda_1 |\innerprod{\hat{\bm x}}{\hat{\bm v}_{\pi^{-1}(1)}} + \innerprod{\hat{\bm x}}{\bm v_1 - \hat{\bm v}_{\pi^{-1}(1)}}| + \eps \nonumber \\
&\le \lambda_1|\innerprod{\hat{\bm x}}{\hat{\bm v}_{\pi^{-1}(1)}}| + \lambda_1 \norm{\bm v_1 - \hat{\bm v}_{\pi^{-1}(1)}}{} + \eps \le \lambda_1 \theta + \cancel{\lambda_1} \frac{(6.2 + 4\kappa) \eps}{\cancel{\lambda_1}} + \eps \nonumber \\
& \le \frac{\lambda_1}{2 \kappa} + 6.2 \eps + 4 \kappa \eps + \eps  \le \frac{\lambda_{\min}}{2} + \frac{\lambda_{\min}}{12.5} + 7.2 \eps < \lambda_{\min} - \eps \nonumber \\
& \le \lambda_{k+1}-\eps, \nonumber
\end{flalign}
where we have used the facts that $
\theta \le 1/(2\kappa)$, $\eps \le \theta^2 \lambda_{\min}/12.5 \le \lambda_{\min}/50$, and
\begin{flalign}
4\kappa \eps \le 4 \cdot \frac{\lambda_{\max}}{\cancel{\lambda_{\min}}}\cdot \frac{\theta^2 \cancel{\lambda_{\min}}}{12.5} \le \frac{\cancel{4} \cdot \lambda_{\max}}{12.5\cdot \cancel{4} \kappa^2} \le \frac{\lambda_{\min}}{12.5}.  \nonumber
\end{flalign}

Therefore, in order to satisfy \eqref{eqn:comb_upp_low}, we must have
\begin{flalign}\label{eqn:max_compare}
\max\{\lambda_1 |x_1|^{p-2}, \lambda_{k+1} |x_{k+1}|^{p-2}\} = \lambda_{k+1}|x_{k+1}|^{p-2},
\end{flalign}
which simplifies \eqref{eqn:comb_upp_low} to
\begin{flalign}\label{eqn:combine_upp_low_2}
\lambda_{k+1} - \eps  \le \max_{i \in [n]\backslash [k]} \lambda_i - \eps  \le \hat \lambda \le \lambda_{k+1} |x_{k+1}|^{p-2} + \eps.
\end{flalign}
Based on  \eqref{eqn:combine_upp_low_2}, we have that
\begin{flalign}\label{eqn:intermediate_bd}
&\lambda_{k+1} \ge \max_{i \in [n]\backslash [k]} \lambda_i - 2 \eps, \quad
|\lambda_{k+1} - \hat{\lambda}| \le \eps, \quad \mbox{and} \\
&|x_{k+1}| \ge |x_{k+1}|^{p-2} \ge \frac{\lambda_{k+1} - 2 \eps}{\lambda_{k+1}} = 1 - \frac{2\eps}{\lambda_{k+1}}. \nonumber
\end{flalign}

Thus, we have achieved the eigenvalue perturbation bound \eqref{eqn:lam_bd} promised in the theorem. Next, we will sharpen the eigenvector perturbation bound by exploiting the optimality conditions for problem \eqref{eqn:subproblem}.

The key observation is that, at the point $\hat{\bm x}$, the constraint $|\innerprod{\hat{\bm v}_i}{\hat{\bm x}}| \le \theta$ is not active. To see this, for any $i \in [k]$,
\begin{flalign}
&|\innerprod{\hat{\bm v}_i}{\hat{\bm x}}| \label{eqn:interior}\\
&= |\innerprod{{\bm v}_{\pi(i)}}{\hat{\bm x}} + \innerprod{\hat{\bm v}{}_i-\bm v_{\pi(i)}}{\hat{\bm x}}| \le |x_{\pi(i)}| + \norm{\hat{\bm v}_i-\bm v_{\pi(i)}}{} \le |x_{\pi(i)}| + (6.2 + 4 \kappa)\eps/\lambda_{\min}  \nonumber \\
& \le \sqrt{1- x_{k+1}^2} + (6.2 + 4 \kappa)\theta^2 /12.5
 \le \sqrt{4 \eps /\lambda_{\min}} + (6.2 + 4 \kappa)\theta/12.5 \cdot \theta \nonumber \\
& \le \sqrt{\frac{4 \theta^2}{12.5}} + \left(\frac{3.1}{12.5} + \frac{2}{12.5}\right) \cdot \theta
< \theta, \nonumber
\end{flalign}
where the last line is due to \eqref{eqn:intermediate_bd} and the fact that $\kappa \theta \le 1/2$ and $\eps \le \theta^2 \lambda_{\min}/12.5.$ Therefore, only the equality constraint is active and will be involved in the optimality conditions at the point $\hat{\bm x}$. Consider the Lagrangian function at the point $\hat{\bm x}$,
\[
\mc L(\hat{\bm x}, \lambda) = \wh{\mcb T} \hat{\bm x}^{\tensor p} - \frac{p \lambda}{2}\left(\norm{\hat{\bm x}}{}^2 - 1\right),
\]
where $\lambda \in \reals$ corresponds to the (scaled) Lagrange multiplier for the equality constraint on the norm of $\hat{\bm x}$, which we have squared. Since the linear independent constraint qualification \cite[Section~12.3]{wright1999numerical} can be easily verified, by the first-order optimality conditions (a.k.a. KKT condition), there exists a $\bar \lambda \in \reals$ such that
\[
\frac{1}{p}\left(\nabla \mc L\left(\hat{\bm x}, \bar \lambda\right)\right) = \wh{\mcb T} \hat{\bm x}^{\tensor p-1} - \bar \lambda \hat{\bm x} = 0.
\]
Moreover, as $\norm{\hat{\bm x}}{}=1$, $\bar \lambda = \bar \lambda \innerprod{\hat{\bm x}}{\hat{\bm x}} = \wh{\mcb T} \hat{\bm x}^{\tensor p} = \hat \lambda$. Thus, we have
\begin{flalign}
\hat \lambda \hat{\bm x} = \wh{\mcb T}\hat{\bm x}^{\tensor {p-1}} = \lambda_{k+1} x_{k+1}^{p-1} \bm v_{k+1} + \sum_{i \neq k+1} \lambda_{i} x_{i}^{p-1} \bm v_{i} + \mcb E  \hat{\bm x}^{\tensor {p-1}}. \nonumber
\end{flalign}
Consider the quantity
\begin{flalign}
& \norm{\lambda_{k+1}(\hat{\bm x} - \bm v_{k+1})}{}  \label{eqn:total} \\
= &\norm{(\lambda_{k+1} - \hat \lambda) \hat{\bm x} + (\hat\lambda \hat{\bm x} - \lambda_{k+1} \bm v_{k+1}) }{}  \nonumber \\
= &  \norm{(\lambda_{k+1} - \hat \lambda) \hat{\bm x} + \lambda_{k+1} (x_{k+1}^{p-1}-1) \bm v_{k+1} + \sum_{i \neq k+1} \lambda_{i} x_{i}^{p-1} \bm v_{i} + \mcb E  \hat{\bm x}^{\tensor {p-1}}}{} \nonumber  \\
\le &  |\lambda_{k+1} - \hat \lambda| + \lambda_{k+1}\vert x_{k+1}^{p-1}-1 \vert  + \norm{\sum_{i \neq k+1} \lambda_{i} x_{i}^{p-1} \bm v_{i}}{} + \norm{\mcb E  \hat{\bm x}^{\tensor {p-1}}}{}. \nonumber
\end{flalign}
Thanks to the intermediate result \eqref{eqn:intermediate_bd}, we have
\begin{flalign}\label{eqn:big_term_1}
&|\lambda_{k+1} - \hat \lambda| \le \eps, \quad  \norm{\mcb E  \hat{\bm x}^{\tensor {p-1}}}{} \le \eps, \quad \mbox{and} \\ &\lambda_{k+1}|x_{k+1}^{p-1}-1| = \lambda_{k+1}\left(1-\Abs{x_{k+1}}\cdot\Abs{x_{k+1}}^{p-2}\right)\le \lambda_{k+1}\left(1- \left(1-2\frac{\eps}{\lambda_{k+1}}\right)^2  \right)\le 4\eps. \nonumber
\end{flalign}
Moreover, for the term $\norm{\sum_{i\neq k+1} \lambda_{i} x_{i}^{p-1} \bm v_{i}}{}$, we can derive that
\begin{flalign}
&\norm{\sum_{i\neq k+1} \lambda_{i} x_{i}^{p-1} \bm v_{i}}{} =  \left( \sum_{i\neq k+1} \lambda_i^2 x_i^{2p-2} \right)^{1/2} \nonumber\\
&\qquad \qquad \le \max \set{\lambda_1 \Abs{x_1}^{p-2}, \lambda_{k+2} \Abs{x_{k+2}}^{p-2}} \sqrt{\sum_{i\neq k+1} x_i^2} \le 4 \kappa\eps。 \label{eqn:big_term_2}
\end{flalign}
The last line holds due to \eqref{eqn:ordering_ass} and for $j \in \set{1, k+2}$,
\begin{flalign}
\lambda_j \Abs{x_j}^{p-2}  \sqrt{\sum_{i\neq k+1} x_i^2} &\le \lambda_j \sqrt{1-x_{k+1}^2} \cdot \sqrt{1-x_{k+1}^2} = \lambda_j (1-x_{k+1}^2) \le \frac{4\lambda_j \eps}{\lambda_{k+1}} \le 4 \kappa \eps,
\end{flalign}
where we have used $\sum_{i \in [n]} x_i^2 = 1$ and \eqref{eqn:intermediate_bd}.

Therefore, by substituting \eqref{eqn:big_term_1} and \eqref{eqn:big_term_2} into \eqref{eqn:total},  one has
\[
\norm{\lambda_{k+1}(\hat{\bm x} - \bm v_{k+1})}{} \le (6 + 4 \kappa) \eps,
\]
which leads to the desired bound
$\norm{\hat{\bm x} - \bm v_{k+1}}{} \le (6.2 + 4 \kappa) \eps/\lambda_{k+1}.$

By mathematical induction, we complete the proof.
\end{proof}

\subsection{Numerical Experiments} \label{subsec:numeric}
In this subsection, we present three sets of numerical experiments to corroborate our theoretical findings in Theorem \ref{sec:SROAwCD} regarding the SROAwCD method. We solve the main subproblem \eqref{eqn:SROAwCD} via the general polynomial solver \textsf{GloptiPoly~3}~\citep{henrion2009gloptipoly}, which is a global solver based on the Sum-of-Squares (SOS) framework \cite{shor87poly,nesterov00squared,parrilo2000structured,lasserre2001global,parrilo2003semidefinite}.

\begin{figure}
\centerline{
\begin{minipage}{3in}
\centerline{\includegraphics[width=1\textwidth]{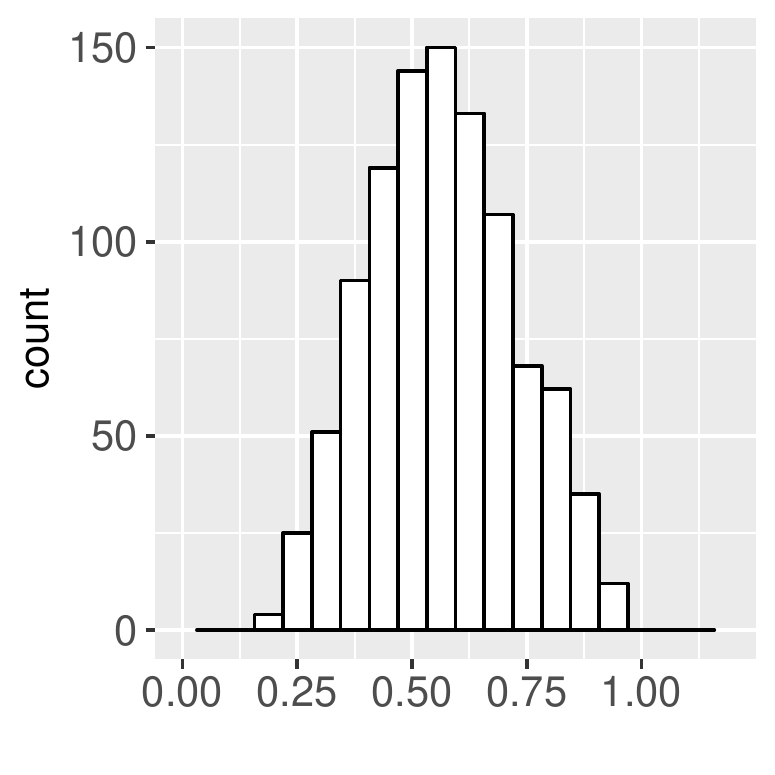}}
\vspace{-2mm}
\centerline{\hspace{20mm}$\frac{1}{\eps} \cdot \max_{j \in [5]} \min \set{\Abs{\hat{\lambda}_j \pm 300}}$}
\end{minipage}
\begin{minipage}{3in}
\centerline{\includegraphics[width=1\textwidth]{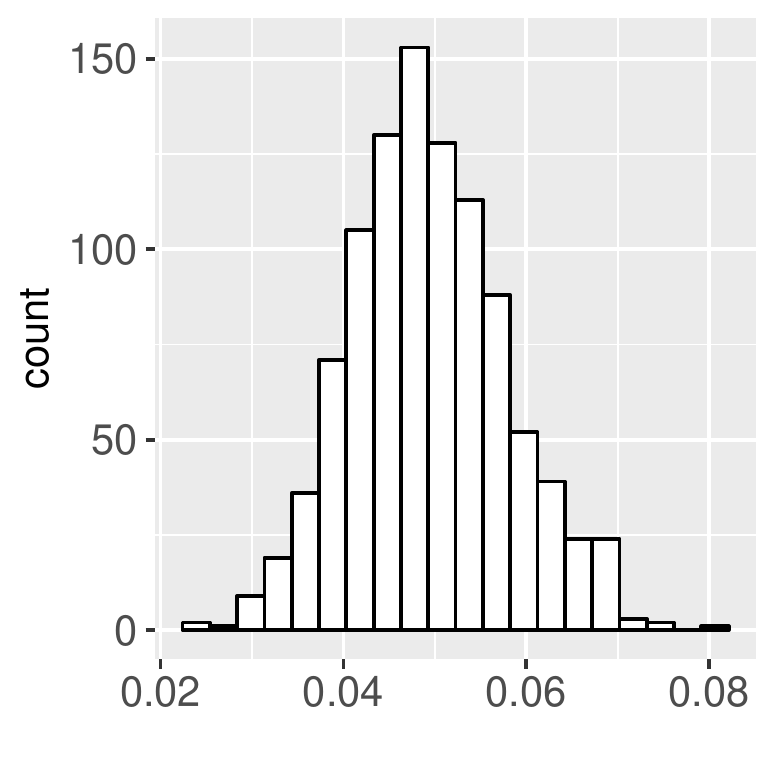}}
\vspace{-2mm}
\centerline{\hspace{15mm} $\frac{300}{10.2 \eps} \cdot \max_{j \in [5]} \min_{i \in [5]} \set{\norm{\hat{\bm v}_j \pm \bm e_i}{}} $}
\end{minipage}
}
\caption{{\bf Histograms for the eigenvalue and eigenvector perturbations in the first experiment.} The left figure plots the histogram of the (normalized) eigenvalue perturbations. All perturbations in the eigenvalues are upper bounded by 1, which is consistent with Theorem \ref{sec:SROAwCD}.  The right figure plots the histogram of the (normalized) eigenvector perturbations. All perturbations in the eigenvectors are upper bounded by 1, which is also consistent with Theorem \ref{sec:SROAwCD}}
\label{fig:1st_exp}
\end{figure}

\paragraph{Experiment 1}
In this experiment, we will synthetically verify the perturbation bounds stated in Theorem \ref{sec:SROAwCD}. We generate nearly symmetric orthogonally decomposable tensor $\widehat{\mcb T} = \mcb T + \mcb E \in \reals^{5 \times 5 \times 5}$ in the following manner. The underlying symmetric orthogonally decomposable tensor $\mcb T$ is set as the diagonal tensor with all diagonal entries equal to 300, i.e., $\mcb T = \sum_{i=1}^5 300 \cdot \bm e_i^{\otimes 3}$, and the perturbation tensors $\mcb E$ are produced by symmetrizing a randomly generated $5 \times 5 \times 5$ tensor whose entries follow standard normal distribution independently.
We set $\theta$ to be $1/(2 \kappa)= 1/2$ (as suggested in Theorem \ref{thm:C-SROA}). 1000 random instances are tested. Figure \ref{fig:1st_exp} plots the histogram of perturbations in both eigenvalue and eigenvector. As depicted in Figure \ref{fig:1st_exp}, both types of perturbations are well controlled by the bounds provided in Theorem \ref{sec:SROAwCD}.

\paragraph{Experiment 2} In Theorem \ref{sec:SROAwCD}, the parameter $\theta$ is suggested to be set to $1/(2\kappa)$. In this experiment, we compare the performance of SROAwCD with $\theta = 1/(2 \kappa) =1/2$ and $\theta = 0$ based on the criterion
\begin{flalign}\label{eqn:metric}
\norm{\mcb T - \sum_{i=1}^5 \hat{\lambda}_i \hat{\bm v}_i^{\otimes 3}}{F}.
\end{flalign}

 The tensors are generated in the same way as in the first experiment. Among all the 1000 random cases, the SROAwCD method with $\theta = 1/2$ consistently outperforms the one with $\theta = 0$. This makes intuitive sense. As $(\hat \lambda_k, \hat{\bm v}_k)$ only approximate the underlying truth, setting $\theta = 0$, which forces strict orthogonality, tends to introduce additional errors into the problem.

\begin{figure}
\centerline{
\begin{minipage}{3in}
\centerline{\includegraphics[width=1\textwidth]{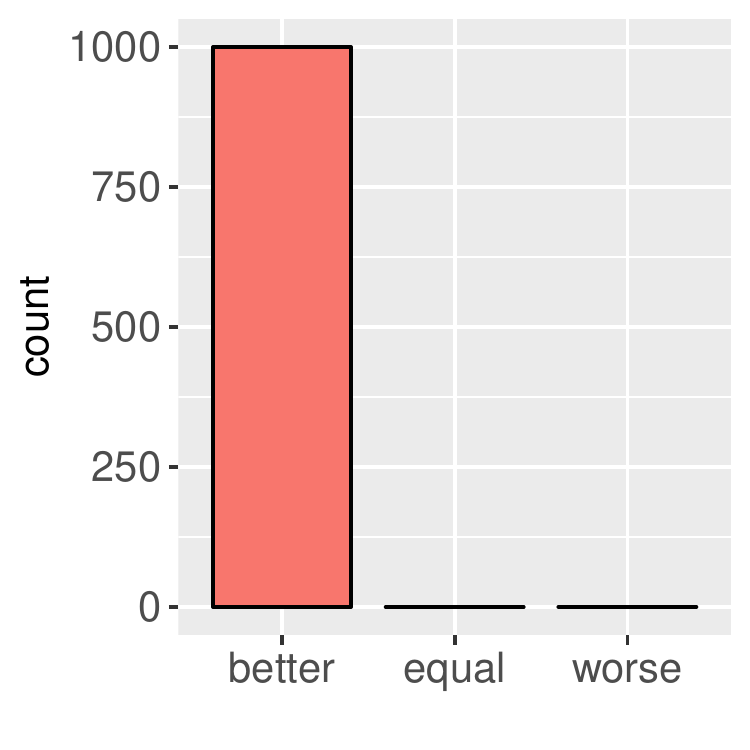}}
\end{minipage}
\begin{minipage}{3in}
\centerline{\includegraphics[width=1\textwidth]{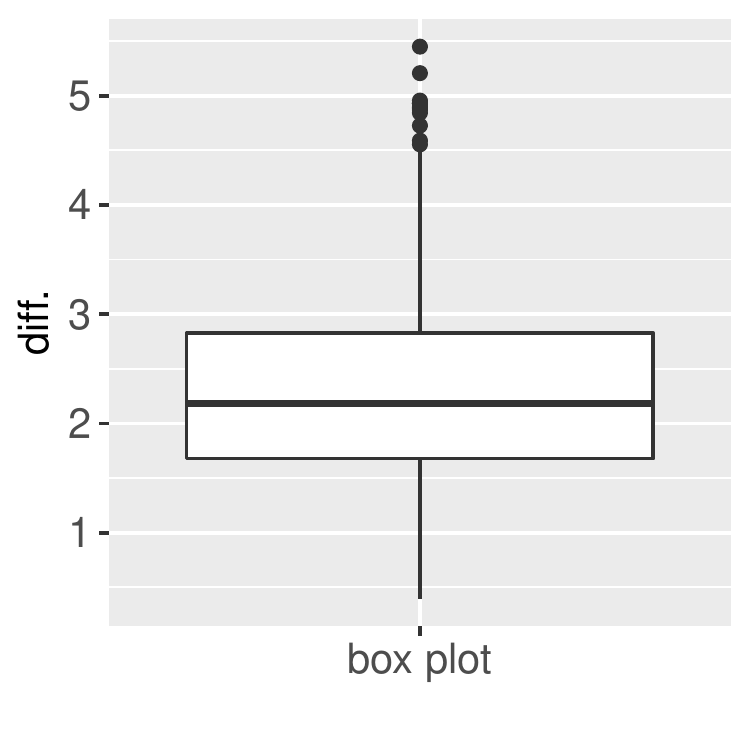}}
\end{minipage}
}
\caption{{\bf Performance comparison for the SROAwCD method with $\theta = 1/2$ and $\theta = 0$ in the second experiment.} As plotted in the left figure, with respect to the criterion \eqref{eqn:metric}, the SROAwCD method with $\theta = 1/2$ outperforms the one with $\theta = 0$ for all randomly generated cases. The right figure is the boxplot of the differences in \eqref{eqn:metric} between the two approaches.}
\end{figure}

\paragraph{Experiment 3} In Theorem \ref{sec:SROAwCD}, the parameter $\theta$ is suggested to be set to $1/(2\kappa)$, which depends on $\kappa$. In this experiment, we will demonstrate the necessity of this dependency. We consider the SOD tensor $\mcb T = 1000 \cdot \bm e_1^{\otimes 3} + \sum_{i = 2}^5 100 \cdot \bm e_i^{\otimes 3}$ with $\mcb E = \mb 0$. We first apply the SROAwCD method with $\theta = 1/2$. The output is as follows:
\begin{flalign*}
&\hat{\bm \lambda} = (1000.00, \;189.95, \; 189.95, \;189.95, \; 189.95)^\top, \quad \mbox{and} \\
&\hat{\bm v}_1 = (1.00, \; 0.00, \; 0.00, \;0.00, \; 0.00)^\top\\
&\hat{\bm v}_2 = (0.50, \; 0.00, \; 0.87, \;0.00, \; 0.00)^\top\\
&\hat{\bm v}_3 = (0.50, \; 0.00, \; 0.00, \;0.87, \; 0.00)^\top\\
&\hat{\bm v}_4 = (0.50, \; 0.87, \; 0.00, \;0.00, \; 0.00)^\top\\
&\hat{\bm v}_5 = (0.50, \; 0.00, \; 0.00, \;0.00, \; 0.87)^\top,
\end{flalign*}
which {\em deviate greatly} from the underlying eigenvalues and eigenvectors of $\mcb T$. Next, we apply the SROAwCD method again with $\theta = {1}/{(2\kappa)} = 1/20$ and the output is as follows:
\begin{flalign*}
&\hat{\bm \lambda} = (1000.00, \;100.00, \; 100.00, \;100.00, \; 100.00)^\top, \quad \mbox{and} \\
&\hat{\bm v}_1 = (1.00, \; 0.00, \; 0.00, \;0.00, \; 0.00)^\top\\
&\hat{\bm v}_2 = (0.00, \; 0.00, \; 1.00, \;0.00, \; 0.00)^\top\\
&\hat{\bm v}_3 = (0.00, \; 1.00, \; 0.00, \;0.00, \; 0.00)^\top\\
&\hat{\bm v}_4 = (0.00, \; 0.00, \; 0.00, \;0.00, \; 1.00)^\top\\
&\hat{\bm v}_5 = (0.00, \; 0.00, \; 0.00, \;1.00, \; 0.00)^\top,
\end{flalign*}
which {\em exactly} recovers (up to a permutation) the underlying eigenvalues and eigenvectors of $\mcb T$.

\subsection{Determination of the maximum spectral ratio $\kappa$}\label{subsec:disc}

As suggested by Theorem \ref{thm:C-SROA}, to choose a proper $\theta$ for the SROAwCD method, we need a rough estimate for $\kappa$. This is not much of a problem, especially for applications in statistics and machine learning, due to several reasons. First, in many problems, we know the maximum spectral ratio $\kappa$ in advance. For example, if we apply independent component analysis \cite{comon1994independent, comon2010handbook} to the dictionary learning model considered in \cite{spielman2012exact, sun2015complete}, it is known that $\kappa = 1$. Moreover, we can always pick the most favorable estimates for $\kappa$ using cross validation \cite[Section~7.10]{friedman2001elements} based on the prediction errors. Furthermore, as a supplement, we can modify Algorithm \ref{alg:SROAwCD} to allow the algorithm to determine the appropriate $\theta$ at each step through adaptive learning. We present the analysis of this modification in the appendix.

\section{Conclusion}
In this paper, we are concerned with finding the (approximate) symmetric and orthogonal decomposition of a nearly symmetric and decomposable (SOD) tensor. Two natural incremental rank-one approximation approaches, the SROAwRD and the SROAwCD methods, have been considered. We first reviewed the existing perturbation bounds for the SROAwRD method. Then we established the first perturbation results for the SROAwCD method that have been given, and discussed issues and potential advantages of this approach. Numerical results were also presented to corroborate our theoretical findings.

We hope our discussion can also shed light on the numerical side. In the SROAwRD method, the main computational bottleneck is the tensor best rank-one approximation problem \eqref{eqn:TBRO}, to which a large amount of attention from a numerical optimization point of view has been paid and for which many efficient numerical methods (e.g. ~\cite{de2000best, zhang2001rank,
kofidis2002best, wang2007successive, kolda2011shifted,
han2012unconstrained, chen2012maximum, zhang2012best, l2014sequential, jiang2012tensor, nie2013semidefinite,
yang2014properties}) have been successfully proposed. In the SROAwCD method, the main computational concern is problem \eqref{eqn:SROAwCD}, which is similar to but slightly more complicated than problem \eqref{eqn:TBRO} with additional linear inequalities. Though general-purpose polynomial solvers based on the sum-of-squares framework \cite{shor87poly, nesterov00squared, parrilo2000structured, lasserre2001global, parrilo2003semidefinite, henrion2009gloptipoly, sostools} can be utilized (as we did in the subsection \ref{subsec:numeric}), more efficient and scalable methods (e.g. projected gradient method \cite{wright1999numerical}, semidefinite programming relaxations \cite{jiang2012tensor, nie2013semidefinite, hu2016note}), specifically tailored to the structure of \eqref{eqn:SROAwCD}, may be anticipated. This is definitely a promising future research direction.

\section*{Acknowledgements}
We are grateful to the associate editor Tamara G. Kolda and two anonymous reviewers for their helpful suggestions and comments that substantially improved the paper.

\appendix

\section{Adaptive SROAwCD}\label{sec:FMSROA}
In this section, we provide a modification to the SROAwCD method that adaptively learns an appropriate $\theta$ at each step based on the information collected so far. The complete algorithm is described in Algorithm \ref{alg:FM-rank-1}.  Note that we initially set $\theta$ as $1/2$ (which is the largest value allowed in Theorem \ref{thm:C-SROA}), and gradually reduce it by checking certain conditions.

\begin{algorithm}
\caption{\underline{Ada}ptive \underline{S}uccessive \underline{R}ank-\underline{O}ne \underline{A}pproximation \underline{w}ith \underline{C}onstrained \underline{D}eflation
(AdaSROAwCD)}
\label{alg:FM-rank-1}
\begin{algorithmic}[1]
  \renewcommand\algorithmicrequire{\textbf{input}}

  \REQUIRE symmetric tensor $\wh{\mcb T} \in \bigotimes^p \reals^n$.
  \STATE \textbf{initialize} $\theta \gets 1/2$ and $\bm v_0 \gets \mb 0$
  \FOR{$k=1$ to $n$}
  \STATE Solve the following optimization problem:
           \begin{flalign}
           (\hat\lambda_k,\hat{\bm v}_k) \in &\argmin_{\lambda\in \reals, \bm v \in \reals^n} \quad \norm{\wh{\mcb T} - \lambda \bm v^{\tensor p}}{F} \label{eqn:C-SROA-FM}\\
           &\;\;\qquad\mbox{s.t.} \quad \norm{\bm v}{} = 1    \nonumber\\
           & \qquad \qquad \;\;\;\; |\innerprod{\bm v}{\hat{\bm v}_i}|\le \theta, \quad i \in [k-1]  \nonumber
           \end{flalign}
           to obtain $(\hat\lambda_k,\hat{\bm v}_k)$.
  \STATE \WHILE{there exists one $i\in [k-1]$ such that $|\innerprod{\hat{\bm v}_k}{\hat{\bm v}_i}| \ge \min\{\hat{\lambda}_k/1.35\hat{\lambda}_i, \theta\}$ } \label{eqn:while_cond}
  \STATE $\theta \gets  0.96\cdot \theta$
  \STATE Solve problem \eqref{eqn:C-SROA-FM} again to replace $(\hat{\lambda}_k, \hat{\bm v}_k)$
  \ENDWHILE
  \STATE $\theta_k \gets \theta$
  \ENDFOR
  \RETURN $\{ (\theta_k, \hat\lambda_k, \hat{\bm v}_k) \}_{k \in [n]}$.
\end{algorithmic}
\end{algorithm}

\begin{theorem}\label{thm:AC-SROA}
Let $\wh{\mcb T} := \mcb T + \mcb E \in \Tensor^p \reals^n$, where $\mcb T$ is a symmetric tensor with orthogonal decomposition $\mcb T
  = \sum_{i=1}^n \lambda_i \bm v_i^{\tensor p}$, $\set{\bm v_1, \bm
  v_2, \dotsc, \bm v_n}$ is an orthonormal basis of $\reals^n$,
  $\lambda_i \neq 0$ for all $i \in [n]$, and $\mcb E$ is a symmetric
  tensor with operator norm $\eps:=\norm{\mcb E}{}$.
  Assume $\eps\le \lambda_{\min}/70\kappa^2$, where  $\kappa := \lambda_{\max}/\lambda_{\min}$, $\lambda_{\min} := \min_{i\in[n]} \Abs{\lambda_i}$ and $\lambda_{\max} := \max_{i\in[n]} \Abs{\lambda_i}$. Then Algorithm \ref{alg:FM-rank-1} terminates in a finite number of steps and its output $\{ (\theta_i, \hat\lambda_i,\hat{\bm v}_i) \}_{i \in [n]}$ satisfies:
  \[
  1/2 = \theta_1 \ge \theta_2 \ge \cdots \ge \theta_n > 0.96/2\kappa,
  \]
  and there exists a permutation $\pi$ of $[n]$ such that for each $j \in [n]$
\begin{flalign*}
&\min\; \set{|\lambda_{\pi(j)} - \hat{\lambda}_j|, \;|\lambda_{\pi(j)} + \hat{\lambda}_j|} \le \eps, \\
&\min\; \set{\norm{\bm v_{\pi(j)} - \hat{\bm v}_j}{},\; \norm{\bm v_{\pi(j)} + \hat{\bm v}_j}{}} \le (6.2 + 4\kappa) \eps/ |\lambda_{\pi(j)}|.
\end{flalign*}
\end{theorem}

\begin{remark}
The condition in line $\ref{eqn:while_cond}$ of Algorithm \ref{alg:FM-rank-1} consists of two components, which are desired properties mainly inspired by the proof of Theorem \ref{thm:C-SROA}.
The first desired inequality $\Abs{\innerprod{\hat{\bm v}_k}{\hat{\bm v}_i}} < \hat{\lambda}_k/1.35\hat{\lambda}_i$ would allow us to establish properties similar to \eqref{eqn:max_compare}. The second desired inequality $\Abs{\innerprod{\hat{\bm v}_k}{\hat{\bm v}_i}} < \theta$ would help us make use of the optimality condition to sharpen the perturbation bounds for the eigenvectors.
The constants chosen in Algorithm \ref{alg:FM-rank-1} and Theorem \ref{thm:AC-SROA}, mainly for illustrative purposes, might be better optimized.
\end{remark}

\begin{proof}
Similar to the proof for Theorem \ref{thm:C-SROA}, without loss of generality, we assume $p\ge 3$ is odd, and $\lambda_i > 0$ for all $i \in [n]$.
Then problem \eqref{eqn:C-SROA-FM} can be equivalently written as
\begin{flalign}\label{eqn:equiv_odd_F}
\hat{\bm v}_k \in \argmax_{\bm v \in \reals^n} \;\wh{\mcb T} \bm v^{\tensor p} \qquad \mbox{s.t.} \quad \norm{\bm v}{} = 1, \; \mbox{and} \; |\innerprod{\bm v}{\hat{\bm v}_i}| \le \theta \;\; \forall \; i \in [k-1],
\end{flalign}
and $\hat \lambda_k = \wh{\mcb T} \hat{\bm v}_k^{\tensor p}$.

Our proof is by induction.

The base case regarding $(\theta_1, \hat{\lambda}_1, \hat{\bm v}_1)$ is the same as the base case in the proof of Theorem \ref{thm:C-SROA}.

We now make the induction hypothesis that for some $k\in [n-1]$, $\{ (\theta_i, \hat\lambda_i,\hat{\bm v}_i) \}_{i \in [k]}$ satisfies
  \[
  1/2 = \theta_1 \ge \theta_2 \ge \cdots \ge \theta_k > 0.96/2\kappa,
  \]
  and there exists a permutation $\pi$ of $[n]$ such that
  \[
    |\lambda_{\pi(j)} - \hat{\lambda}_j| \le \eps, \qquad  \norm{\bm
    v_{\pi(j)}-\hat{\bm v}_j}{} \le (6.2 + 4 \kappa) \eps/\lambda_{\pi(j)}, \quad
    \forall j \in [k].
  \]
Then we are left to prove that
\begin{flalign}
0.96/2 \kappa <\theta_{k+1}\le \theta_k,
\end{flalign}
and there exists an $l \in [n] \backslash \set{\pi(j): j\in [k]}$ that satisfies
\begin{flalign}\label{eqn:indhyp_FM}
|\hat \lambda_{k+1} - {\lambda}_l|\le \eps, \quad \mbox{and} \quad \norm{\hat{\bm v}_{k+1}- {\bm v}_l}{}
      \le \frac{(6.2 + 4\kappa) \eps}{\lambda_{\pi(l)}}.
\end{flalign}

To prove that $\theta_{k+1} > \frac{0.96}{2 \kappa}$, we show that whenever $\theta \in (\frac{0.96}{2\kappa}, \frac{1}{2\kappa}]$, the condition for the while loop in line $\ref{eqn:while_cond}$ of Algorithm \ref{alg:FM-rank-1} will always be satisfied, so $\theta$ can never be reduced to any value below $\frac{0.96}{2\kappa}$.

Consider $\theta \in (\frac{0.96}{2\kappa}, \frac{1}{2\kappa}]$, and denote
\begin{flalign}\label{eqn:rank-one-FM}
\hat{\bm x} \in \argmax_{\bm v \in \reals^n} \;\wh{\mcb T} \bm v^{\tensor p} \qquad \mbox{s.t.} \quad \norm{\bm v}{} = 1, \; \mbox{and} \; |\innerprod{\bm v}{\hat{\bm v}_i}| \le \theta \;\; \forall \; i \in [k],
\end{flalign}
and $\hat \lambda = \wh{\mcb T} \hat{\bm x}^{\tensor p}$.
As $\theta \le 1/2\kappa$ and $\eps \le \lambda_{\min}/70\kappa^2 \le \frac{0.96^2 \lambda_{\min}}{4\cdot 12.5 \cdot \kappa^2} \le \theta^2 \lambda_{\min}/12.5$, results from Theorem \ref{thm:C-SROA} and its proof can be directly borrowed. Based on \eqref{eqn:v_feas}, we know that for any $i \in \set{\pi(j) \; \vert \; j \in [n]\backslash[k]}$, $\bm v_i$ is feasible to problem \eqref{eqn:rank-one-FM}. Then it can be easily verified that
\begin{flalign}\label{eqn:lower_bd_1}
\hat \lambda \ge  \max_{i \in \set{\pi(j) \; \vert \; j \in [n]\backslash[k]}} \; (\mcb T + \mcb E)\bm v_i^{\tensor p} \ge \lambda_{\min} - \eps.
\end{flalign}
Moreover, for any $i\in [k]$,
\begin{flalign}\label{eqn:upper_bd_1}
\hat \lambda_i \;\le\; \max_{\norm{\bm v}{}=1}\; \mcb  T \bm v^{\tensor p}
 + \max_{\norm{\bm v}{}=1}\; \mcb  E \bm v^{\tensor p} = \lambda_{\max} + \eps. \end{flalign}
Using \eqref{eqn:lower_bd_1} and \eqref{eqn:upper_bd_1}, we obtain that for each $i \in [k]$,
\begin{flalign}\label{eqn:condition_1}
\frac{\hat \lambda}{1.35 \hat \lambda_i} \ge \frac{\lambda_{\min}- \eps}{1.35(\lambda_{\max} + \eps)} \ge \frac{\lambda_{\min} - \lambda_{\min}/70\kappa^2}{1.35(\lambda_{\max} + \lambda_{\min}/70\kappa^2)} \ge \frac{69/70\cdot \lambda_{\min}}{1.35\cdot 71/70 \cdot \lambda_{\max}} > \frac{1}{2\kappa} \ge \theta
\end{flalign}
Based on \eqref{eqn:interior}, we also know that for all $i\in [k]$,
\begin{flalign}\label{eqn:condition_2}
|\innerprod{\hat{\bm v}_i}{\hat{\bm x}}| < \theta.
\end{flalign}
So, the combination of \eqref{eqn:condition_1} and \eqref{eqn:condition_2} leads to
\begin{flalign}
|\innerprod{\hat{\bm v}_i}{\hat{\bm x}}| < \theta = \min \set{\theta, \frac{\hat \lambda}{1.35 \hat \lambda_i} }, \quad \forall i\in [k],
\end{flalign}
which implies that $\hat{\bm x}$ satisfies the condition in the while loop. Therefore, as we argued previously, we must have $\theta_{k+1}> \frac{0.96}{2\kappa}$. So $\theta_{k+1}$ is either in $(\frac{0.96}{2\kappa}, \frac{1}{2\kappa}]$ or in $(\frac{1}{2\kappa}, \frac{1}{2}]$.

For the first case, i.e. $\theta_{k+1} \in (\frac{0.96}{2\kappa}, \frac{1}{2\kappa}]$, we can directly establish the result by using the argument for the induction hypothesis \eqref{eqn:indhyp} in the proof of Theorem \ref{thm:C-SROA}.

Hence in the following, we only focus on the second case where $\theta_{k+1} \in (\frac{1}{2\kappa}, \frac{1}{2}]$.

Denote $\hat{\bm x} = \sum_{i \in [n]} x_i \bm v_i:= \hat{\bm v}_{k+1}$ and $\hat \lambda := \hat \lambda_{k+1}.$ Without loss of generality, we renumber $\set{\left(\lambda_{\pi(i)}, \bm v_{\pi(i} \right)}_{i \in [k]}$ to $\set{\left(\lambda_i, \bm v_i \right)}_{i \in [k]}$ and renumber $\set{\left(\lambda_i, \bm v_i \right)}_{i \in [n] \backslash \set{\pi(i) \vert i \in [k]}}$ to $\set{\left(\lambda_i, \bm v_i \right)}_{i \in [n]\backslash[k]}$, respectively, to satisfy
\begin{flalign}\label{eqn:ordering_ass_FM}
&\lambda_1 |x_1|^{p-2} \ge \lambda_2 |x_2|^{p-2} \ge \ldots \ge \lambda_k |x_k|^{p-2}, \quad \mbox{and}  \\
&\lambda_{k+1}|x_{k+1}|^{p-2} \ge \lambda_{k+2}|x_{k+2}|^{p-2} \ge \ldots \ge \lambda_{n}|x_{n}|^{p-2}. \nonumber
\end{flalign}

In the following, we will show that $l = k+1$ is the index satisfying \eqref{eqn:indhyp_FM}.

Based on \eqref{eqn:equiv_odd_F},
\begin{flalign}\label{eqn:subproblem_FM}
 \hat{\bm x} \in \arg \min_{ \bm v \in \reals^n} \;\; \wh{\mcb T} {\bm v}^{\tensor p}  \quad \mbox{subject to} \quad \norm{\bm v}{} = 1, \;\;\; |\innerprod{\hat{\bm v}_i}{\bm v}| \le \theta_{k+1} \mbox{ for any }  i \in [k],
\end{flalign}
and $\hat \lambda = \wh{\mcb T} \hat{\bm x}^{\tensor p}.$
We now bound $\hat \lambda$ from below and above.

We first consider the lower bound by finding a $\bm v$ that is feasible for \eqref{eqn:subproblem_FM}. For any $(i,j) \in [n]\backslash [k] \times [k]$, one has
\begin{flalign*}
|\innerprod{\bm v_i}{\hat{\bm v}_j}| &= |\innerprod{\bm v_i}{\bm v_{\pi(j)} + \hat{\bm v}_j - \bm v_{\pi(j)}}| = |\innerprod{\bm v_i}{\bm v_{\pi(j)}} + \innerprod{\bm v_i}{\hat{\bm v}_j - \bm v_{\pi(j)}}| = |\innerprod{\bm v_i}{\hat{\bm v}_j - \bm v_{\pi(j)}}| \\
&\le \norm{\hat{\bm v}_j - \bm v_{\pi(j)}}{} \le (6.2 + 4\kappa) \eps / \lambda_{\pi(j)} \le (6.2 + 4 \kappa) \lambda_{\min} / (70 \kappa^2 \lambda_{\min}) < \frac{1}{2\kappa} \le \theta_k.
\end{flalign*}
Hence, $\set{\bm v_i}_{i \in [n]\backslash [k]}$ are all feasible to problem \eqref{eqn:subproblem_FM} and then we can easily achieve a lower bound for $\hat \lambda$, as
\begin{flalign}\label{eqn:lower_bd_FM}
\hat \lambda = \wh{\mcb T} \hat{\bm x}^{\tensor p} \ge \max_{i \in [n]\backslash [k]} \wh{\mcb T} \bm v_i^{\tensor p} \ge \max_{i \in [n]\backslash [k]} \lambda_i - \eps \ge \lambda_{k+1} - \eps.
\end{flalign}

Regarding the upper bound, one has
\begin{flalign}
\hat \lambda = \wh{\mcb T} \hat{\bm x}^{\tensor p} = (\mcb T + \mcb E) \hat{\bm x}^{\tensor p} & =  \left(\sum_{i=1}^n \lambda_i \bm v_i^{\tensor p} + \mcb E \right) \hat{\bm x}^{\tensor p} \le \max\{\lambda_1 |x_1|^{p-2}, \lambda_{k+1} |x_{k+1}|^{p-2}\} + \eps,  \label{eqn:upper_bd_FM_temp}
\end{flalign}
as in \eqref{eqn:upper_bd}.

Combining \eqref{eqn:lower_bd_FM} and \eqref{eqn:upper_bd_FM_temp}, we have
\begin{flalign}\label{eqn:comb_upp_low_FM}
\lambda_{k+1} - \eps \le \max_{i \in [n]\backslash [k]} \lambda_i - \eps  \le \hat \lambda \le \max\{\lambda_1 |x_1|^{p-2}, \lambda_{k+1} |x_{k+1}|^{p-2}\} + \eps.
\end{flalign}
Also note that
\begin{flalign}
\lambda_1 |x_1|^{p-2} + \eps &\le  \lambda_1 |x_1| + \eps
\le \lambda_1 |\innerprod{\hat{\bm x}}{\bm v_1}| + \eps
= \lambda_1 |\innerprod{\hat{\bm x}}{\hat{\bm v}_{\pi^{-1}(1)}} + \innerprod{\hat{\bm x}}{\bm v_1 - \hat{\bm v}_{\pi^{-1}(1)}}| + \eps \nonumber \\
&\le \lambda_1|\innerprod{\hat{\bm x}}{\hat{\bm v}_{\pi^{-1}(1)}}| + \lambda_1 \norm{\bm v_1 - \hat{\bm v}_{\pi^{-1}(1)}}{} + \eps \le \lambda_1 \frac{\hat \lambda}{1.35 \hat \lambda_{\pi^{-1}(1)}} + (6.2 + 4\kappa) \lambda_1 \eps/\lambda_{1} + \eps < \hat \lambda.  \nonumber
\end{flalign}
Here, we have used $\eps \le \frac{\hat \lambda}{69}$, due to $\eps \le \frac{\lambda_{\min}}{70}$ and $\hat \lambda \ge \lambda_{\min} - \eps.$
Therefore, in order to satisfy \eqref{eqn:comb_upp_low_FM}, we must have
\begin{flalign}\label{eqn:max_compare_FM}
\lambda_{k+1}|x_{k+1}|^{p-2} = \max\{\lambda_1 |x_1|^{p-2}, \lambda_{k+1} |x_{k+1}|^{p-2}\},
\end{flalign}
which simplifies \eqref{eqn:comb_upp_low_FM} to
\begin{flalign}\label{eqn:upper_bd_FM}
\lambda_{k+1} - \eps  \le \max_{i \in [n]\backslash [k]} \lambda_i - \eps  \le \hat \lambda \le \lambda_{k+1} |x_{k+1}|^{p-2} + \eps.
\end{flalign}
Hence,
\begin{flalign*}
&\lambda_{k+1} \ge \max_{i \in [n]\backslash [k]} \lambda_i - 2 \eps, \quad
|\lambda_{k+1} - \hat{\lambda}| \le \eps, \quad \mbox{and} \\
&|x_{k+1}| \ge |x_{k+1}|^{p-2} \ge \frac{\lambda_{k+1} - 2 \eps}{\lambda_{k+1}} = 1 - \frac{2\eps}{\lambda_{k+1}}.
\end{flalign*}

The eigenvector perturbation bound can be sharpened by exploiting the optimality condition of problem \eqref{eqn:subproblem_FM} as what we did in the proof of Theorem \ref{thm:C-SROA}.  As explicitly required in the algorithm, the constraint $|\innerprod{\hat{\bm v}_i}{\hat{\bm x}}| \le \theta_k$ is not active for any $i \in [k]$ at the point $\hat{\bm x}$. Then by the optimality condition, we again have
\begin{flalign}
\hat \lambda \hat{\bm x} = \hat{\mcb T}\hat{\bm x}^{\tensor {p-1}}. \nonumber
\end{flalign}
By applying exactly the same argument as in the proof of Theorem \ref{thm:C-SROA}, we can obtain
\[
\norm{\lambda_{k+1}(\hat{\bm x} - \bm v_{k+1})}{} \le (6.2 + 4 \kappa)\eps,
\]
which leads to
$\norm{\hat{\bm x} - \bm v_{k+1}}{} \le (6.2 + 4 \kappa) \eps/\lambda_{k+1}.$

By mathematical induction, we have completed the proof.
\end{proof}

\bibliographystyle{ieeetr}
\bibliography{rank-1}

\end{document}